# A universal preprocessing algorithm of average kernel method with Gauss-Laguerre quadrature for double integrals

Kejun Pan

Key Laboratory of High-Performance Ship Technology, Ministry of Education, Wuhan University of Technology, Wuhan 430063, China

Mingliang Xie

State Key Laboratory of Coal Combustion, Huazhong University of Science and Technology, Wuhan 430074, China
Corresponding author's Email: mlxie@mail.hust.edu.cn

**Abstract:** To address the computational challenges posed by nonlinear collision kernels in the Smoluchowski equation, this study proposes a universal preprocessing algorithm for the average kernel method based on the Gauss-Laguerre quadrature for double integrals. With this algorithm, the numerical code accurately and efficiently determines the pre-exponential factor of the average kernel. Additionally, the exact pre-exponential factors of the four fundamental average kernels and their associated truncation error estimations were analyzed. The results demonstrate the reasonability and reliability of the preprocessing algorithm.

**Introduction**

The Smoluchowski equation is one of the basic equations in aerosol science and technology and is also one of the main equations of kinetic molecular theory, which has a wide range of applications in the fields of mathematics and engineering. For a mono-variant problem, the classical Smoluchowski equation takes the form (**Friedlander, 2000**)

$$\frac{\partial n(v,t)}{\partial t} = \frac{1}{2}\int_0^v \beta(v_1, v-v_1)n(v_1)n(v-v_1)dv_1 - \int_0^\infty \beta(v_1, v)n(v_1)n(v)dv_1 \quad (1)$$

where $n(v,t)dv$ is the number density of particles per unit spatial volume with particle volume from $v$ to $v + dv$ at time $t$, and $\beta$ is the kernel of coagulation.

In the past 100 years, the Smoluchowski equation has been solved analytically for only a limited number of known simple collision kernels, such as constant, additive, and multiplicative kernels (**Leyvraz, 2003**). In the real world, the collision kernel under different physical conditions typically takes a nonlinear form. The convolution in the Smoluchowski equation is usually impossible to integrate directly if the collision kernel is a function of particle size.

In 1940, Schumann proposed the concept of the average kernel to approximate the actual collision kernel, which is defined with Laplace transformation as (**Schumann, 1940**)

$$\int_0^\infty \int_0^\infty \bar{\beta} \exp\left(-\frac{v_1+v}{u}\right) dv_1 dv = \int_0^\infty \int_0^v \beta(v_1, v) \exp\left(-\frac{v_1+v}{u}\right) dv_1 dv \quad (2)$$

where $u$ is the algebraic mean volume of particle size distribution, and it is defined as

$$u = \frac{M_1}{M_0} \tag{3}$$

where $M_0$ and $M_1$ are the moments of volume-based PSD as

$$\begin{cases} M_0(t) = \int_0^\infty n(v,t)dv \\ M_1(t) = \int_0^\infty vn(v,t)dv \end{cases} \tag{4}$$

Through operation and reorganization, the average kernel can be represented as

$$\bar{\beta} = \frac{1}{u^2} \int_0^\infty \int_0^v \beta(v_1, v) \exp\left(-\frac{v_1+v}{u}\right) dv_1 dv \tag{5}$$

For a homogeneous collision kernel (i.e., $\beta(\alpha v, \alpha v_1) = \alpha^q \beta(v, v_1)$, $\beta(v, v_1) = \beta(v_1, v)$ and $\beta(v, v_1) \geq 0$), if the scale factor is $\alpha = 1/u$, the average kernel can be simplified as

$$\bar{\beta} = pu^q \tag{6}$$

where $q$ is power index, $p$ is a proportional factor representing the total collision frequency of particle coagulating system, and it can be calculated as

$$p = \int_0^\infty \int_0^{u\eta} e^{-\eta-\eta_1} \beta(\eta, \eta_1) d\eta_1 d\eta \tag{7}$$

where the dimensionless particle volume is defined as $\eta = v/u$. Due to the symmetry of the homogeneous collision kernel, $p$ can be simplified as

$$p = \frac{1}{2} \int_0^\infty \int_0^\infty e^{-\eta-\eta_1} \beta(\eta, \eta_1) d\eta_1 d\eta \tag{8}$$

Substituting the average collision kernel $\bar{\beta} = pu^q$ into the Smoluchowski equation, the analytical solution can be obtained as

$$n(v,t) = \frac{M_0^2}{M_1} \exp\left[-\frac{M_0}{M_1}v\right] \tag{9}$$

Generally, the zeroth moment ($M_0$) represents the total particle number concentration, and the first moment ($M_1$) is proportional to the particle mass concentration, which remains constant owing to rigorous mass conservation requirements. Then the evolution of zeroth moment with average kernel can be found as

$$\frac{dM_0}{dt} = -\frac{1}{2}\bar{\beta}M_0^2 \tag{10}$$

According to the value of the index $q$, the solution of the zeroth-order moment can be divided into three types (**Xie, 2024**). The solution of the zeroth-order moment together with Eq. (8) provides a complete solution of the Smoluchowski equation, which significantly simplifies its analytical solution procedure of the Smoluchowski equation. Schumann's work on solving PBE with average kernel method has become a classic in this field. So far, this result has been widely cited as a benchmark by scholars in different fields (**Lin et al., 2018; Makoveeva and Alexandrov, 2023**).

The average kernel method offers several notable advantages that make it an attractive choice for many applications. Its primary strengths lie in its simplicity and ease of implementation. This approach provides a straightforward method to combine information from multivariate kernel functions without the need for extensive parameter tuning. The method is particularly useful when prior knowledge of the relative importance of different kernels is limited, as it treats multivariate kernels equally and can still achieve competitive results. In addition, the average kernel approach maintains the positive definiteness of the resulting kernel matrix, ensuring mathematical consistency and applicability to various kernel-based algorithms. This combination of simplicity, robustness, and mathematical soundness makes it a practical choice for many real-world applications, particularly when computational resources are limited or when a quick, reliable solution is required. The average kernel approach, which is useful in certain scenarios, has several significant limitations that make it inuniversally applicable. First, the method may lose important local structure information through averaging and is not well suited for highly nonlinear or multimodal kernels. The analytical solution of the Smoluchowski equation is independent of the average kernel and boundary conditions. Second, there is an irreconcilable contradiction between the analytical solution and the experimental data. Although Schumann adopted the radius-based particle size distribution instead of the volume-based particle size distribution to achieve matching between the analytical solution and experimental data, this treatment did not solve the problem fundamentally. Third, the binary Laplace transformation for a nonlinear collision kernel has higher computational complexity. This approach also requires careful parameter tuning and lacks clear guidelines for optimal multiple-kernel weight selection. Therefore, the average kernel method can be useful for moderate-sized datasets, where kernels capture related aspects and computational efficiency is critical.

Recently, an AK-iDNS framework was proposed to solve the Smoluchowski equation (**Pan *et al.*, 2024**). This framework uses the analytical solution of the Smoluchowski equation based on the average kernel method (AK) as the initial condition, and the self-preserving size distribution is achieved through an iterative direct numerical simulation algorithm (iDNS) (**Xie and He, 2022; Xie, 2023**). The corrected similarity solution is not only consistent with the experimental data but also consistent with the properties of the original kernel. Therefore, the first two defects of the average kernel method were effectively addressed. Nevertheless, we treated the definition of the average kernel as a Laplace transform in a previous work, where the binary Laplace transform projects a spatial surface onto a point on the symmetry plane. Utilizing the homogeneity of the kernel, the projection point is replaced by an intersection point between the surface and symmetric plane. This treatment greatly simplifies the computation of the average kernel and achieves the same scaling law as that of the Taylor series expansion method of moments (**Xie and Wang, 2013; Xie, 2015; Yu et al., 2008**). However, the simplified computational method of the average kernel may produce deviations to some extent and is inconsistent with the physical meaning defined by the average kernel. For example, the average kernel of sedimentation coagulation is zero, which is unreasonable. According to the structure of the integrals, it is natural to consider approximating it with the Gauss-Laguerre quadrature. However, numerical methods and error estimates specifically designed for double integrals of this form have not been extensively studied in the existing literature. An alternative method for obtaining the population-averaged coagulation kernel is provided in **Appendix A**. Because the dimensionless analytical solution of the Smoluchowski equation with an averaged kernel is $\psi(\eta) = e^{-\eta}$, then the average kernel based on Laplace transformation is equivalent to the population-averaged kernel. Therefore, accurate preprocessing methods of average kernel have important theoretical value and broad application prospects

The purpose of this study is to address and improve the third drawback of the average kernel method and provide an effective numerical method for calculating the pre-exponential factor of the average kernel. In this study, the definition of the average kernel was reformulated using the Gaussian-Laguerre quadrature for double integrals. A MATLAB implementation of the proposed method is provided, offering an efficient solution for calculating the double Gauss-Laguerre integrals and their associated errors.

**The Gauss-Laguerre quadrature and its remainder**

Laguerre orthogonal polynomials have been widely used in various applications (**Fatyanov and Terekhov, 2011**). In this study, only the case in which two-dimensional expansion is preserved is considered.

First, the single integral $I$ with the $(n+1)$-point Gauss-Laguerre quadrature rule is given by

$$I = \int_0^\infty e^{-x} f(x)dx = Q_{n+1} + R_{n+1}(f) \tag{11}$$

and
$$Q_{n+1} = \sum_{j=1}^{n+1} w_i f(x_i) \tag{12}$$

where the abscissas $x_i$ are the zeros of the Laguerre polynomials $L_{n+1}(x)$, and $L_{n+1}(x)$ is defined as

$$L_{n+1}(x) = \frac{e^x}{(n+1)!} \frac{d^{n+1}(x^{n+1}e^{-x})}{dx^{n+1}} \tag{13}$$

and $w_i$ are the corresponding weights or Christoffel numbers, and $w_i$ can be expressed as

$$w_i = \int_{-1}^{1} \prod_{i=0, i \neq k}^{n+1} \frac{x-x_i}{x_k-x_i} e^{-x} dx \tag{14}$$

On the assumption that the definition of $f$ may be continued into a complex $z$-plane, where $z = x + iy$. The remainder (or truncation error) of the quadrature rule is then a contour integral. The contour $C_z$ satisfied $Re(z) \in (0, \infty)$ so that the function $f$ satisfies the condition that the required integral on and within $C_z$ exists in the Riemann sense. The remainder is given by **Donaldson and Elliott (1972)** as follows:

$$R_{n+1}(f) = \frac{1}{2\pi i} \oint k_{n+1}(z) f(z) dz \tag{15}$$

in which the function $k_{n+1}(z)$ is given by:

$$k_{n+1}(z) = \frac{\Pi_{n+1}(z)}{L_{n+1}(z)} \tag{16}$$

where the function $\Pi_{n+1}(z)$ can be expressed as

$$\Pi_{n+1}(z) = -\Gamma(n+2)U(n+2,1;e^{-i\pi}) \tag{17}$$

where $\Gamma$ is the Euler gamma function, and $U$ is the second type of confluent hypergeometric function that satisfies $0 < \arg(z) < 2\pi$.

Second, the double integral $II$ with the $(n+1) \times (n+1)$-point Gauss-Laguerre quadrature rule is given by:

$$II = \int_0^\infty \int_0^\infty e^{-x-x_1} f(x, x_1) \, dx_1 dx$$

$$= \sum_{i=1}^{n+1} \sum_{j=1}^{n+1} w_i w_j f(x_i, x_j) + \int_0^\infty R_{n+1,x_1} \, dx_1 + \int_0^\infty R_{n+1,x} \, dx + O(R_{n+1,x}; R_{n+1,x_1}) \quad (18)$$

where the remainders are given by:

$$R_{n+1,x_1} = \frac{1}{2\pi i} \oint k_{n+1}(z) f(x, z) dz \tag{19a}$$

$$R_{n+1,x} = \frac{1}{2\pi i} \oint k_{n+1}(z) f(z, x_1) dz \tag{19b}$$

$$O(R_{n+1,x}; R_{n+1,x_1}) = -\frac{1}{2\pi i} \oint k_{n+1}(z) R_{n+1,x_1}(z) dz - \frac{1}{2\pi i} \oint k_{n+1}(z) R_{n+1,x}(z) dz \tag{19c}$$

Thus, it is assumed that the remainder $R_{n+1,x_1}$ may continue into a complex $z$-plane. In any case, the last term $O(R_{n+1,x}; R_{n+1,x_1})$ is essentially the remainder of the remainder. Therefore, it was assumed to be negligible and was replaced by zero (**Elliott et al., 2011**). On combing these results, the double integral can be expressed as

$$II = \int_0^\infty \int_0^\infty e^{-x-x_1} f(x, x_1) dx_1 dx = Q_{n+1,n+1} + R_{n+1,n+1} \tag{20}$$

in which $Q_{n+1,n+1}$ and $R_{n+1,n+1}$ are given as follows:
$$Q_{n+1,n+1} = \sum_{i=1}^{n+1} \sum_{j=1}^{n+1} w_j w_i f(x_i, x_j) \tag{21}$$
and

$$R_{n+1,n+1} = 2 \int_0^\infty R_{n+1}(f(x,x)) dx + O(R_{n+1,x}; R_{n+1,x_1}) \tag{22}$$

**Simplified error analysis**

In the field of convergence acceleration of sequences and error estimation, the application of the $\theta$-procedure to the rules of a sequence transformation leads to a new sequence transformation (**Brezinski, 1982; Fdil, 1997**). Donaldson and Elliott derived a priori error estimates using contour integration methods (**Donaldson and Elliott, 1972**). Similarly, **Davis and Rabinowitz** (1984) established asymptotic convergence bounds for Laguerre expansions. **Kzaz** (1999) developed asymptotic forms of the remainder estimates for specific function classes. Recent studies (**Xiang, 2012; Zhang et al., 2024**) have further refined the error estimates in quadrature formulas and examined pointwise error bounds for Laguerre polynomial expansions. Building on these foundational studies, this work identifies a comparable power-law relationship in the data through a logarithmic transformation and extends the analysis to a two-dimensional Gaussian-Laguerre quadrature.

If the function $f$ verifies, uniformly in $\arg z$ for $Re(z) \to +\infty$, the condition

$$f(z) = o\left(\frac{e^z}{z^m}\right), \forall m > 0 \tag{23}$$

The truncation error of single Gauss-Laguerre quadrature can be approximated as (**Brezinski, 1982; Fdil, 1997; Kzaz, 1999**)

$$R_{n+1}(f) \sim \frac{1}{\ln^a n} \left(\alpha_0 + \sum \frac{\alpha_{r,s}}{n^r \ln^s n}\right), (a > 0, \ a \in R, \alpha_0 \neq 0) \tag{24}$$

where $\alpha_0$ is a constant related to the number of interpolation nodes $(n)$, $r$ and $s$ are the order parameters. To analyze the truncation error of single Gauss-Laguerre quadrature, a new error variable is introduced here, which is defined as

$$\epsilon_n = |Q_{n+1} - Q_n| \qquad (25)$$

If the number of Laguerre polynomial nodes is large enough, and the error can be expressed approximately as

$$\epsilon_n \sim \frac{\alpha_0}{\ln^a n} \qquad (26)$$

This implies that in the logarithmic coordinate system, the error $\epsilon_n$ is approximately linearly related to the number of interpolation nodes, and decreases with an increase in the number of nodes.

Analogously, the error of double Gauss-Laguerre quadrature is defined as

$$\epsilon_{n,n} = |Q_{n+1,n+1} - Q_{n,n}| \qquad (27)$$

And it can be calculated approximately as

$$\epsilon_{n,n} \sim \frac{\alpha_0}{\ln^a n} \epsilon_n \qquad (28)$$

This means that in the logarithmic coordinate system, the error $\epsilon_{n,n}$ is also approximately linearly related to the number of interpolation nodes, and decreases with an increase in the number of nodes. Without loss of generality, according to the two-point formula of the line, it means that the following relationship holds

$$\frac{\ln \epsilon_{m,m} - \ln \epsilon_{n,n}}{\ln m - \ln n} = C, \quad (m, n \in N) \qquad (29)$$

where $C$ is the constant slope. And it can be reorganized into

$$\frac{\epsilon_{m,m}}{\epsilon_{n,n}} = \left(\frac{m}{n}\right)^C \qquad (30)$$

It can be extended to the real number field as

$$\epsilon_{x,x} = \epsilon_{n,n} \left(\frac{x}{n}\right)^C, \quad (x \in R) \qquad (31)$$

And the truncation error of double Gauss-Laguerre quadrature can be calculated as

$$R_{n+1,n+1} = \int_{n+1}^{\infty} \epsilon_{x,x} dx = -\epsilon_{n,n} \left(\frac{n+1}{n}\right)^C \left(\frac{n+1}{C+1}\right); \quad (C < -1, C \in R) \qquad (32)$$

It should be pointed out that the problem of orthogonal stability of Laguerre polynomials will arise (**Kahaner and Monegato, 1978; Wünsche, 2000**) as the number of interpolation nodes increases, and the zeros of Laguerre polynomials together with positive weights will enter the complex domain at higher-order terms. Therefore, the maximum number of interpolation points is used as $n = 360$ in the present study.

**Results and Discussions**

In this section, several examples are presented to demonstrate the accuracy of the Gauss-Laguerre quadrature for a double integral. Some common collision kernels for various mechanisms of particle interactions are listed bellows,

$$\beta_{CR} = \left(\eta_1^{-\frac{1}{3}} + \eta^{-\frac{1}{3}}\right)\left(\eta_1^{\frac{1}{3}} + \eta^{\frac{1}{3}}\right) \qquad (33)$$

for Brownian coagulation in the continuum regime (CR),

$$\beta_{SC} = \left(\eta_1^{\frac{1}{3}} + \eta^{\frac{1}{3}}\right)^3 \tag{34}$$

for shear coagulation (SC),

$$\beta_{FM} = (\eta_1^{-1} + \eta^{-1})^{\frac{1}{2}} \left(\eta_1^{\frac{1}{3}} + \eta^{\frac{1}{3}}\right)^2 \tag{34}$$

for Brownian coagulation in the free molecule regime (FM),

$$\beta_{SD} = \left(\eta^{\frac{1}{3}} + \eta_1^{\frac{1}{3}}\right)^3 \left|\eta^{\frac{1}{3}} - \eta_1^{\frac{1}{3}}\right| \tag{35}$$

for sedimentary coagulation under gravity (SD), respectively. For brevity, some physical constants were ignored in the expression of the collision kernel. Obviously, all these kernels satisfy the conditions of Eq. (23). The exact expressions of the double integrals and corresponding numerical data based on the average kernel method for the above four types of kernels are listed below.

**Brownian coagulation in the continuum regime**

The double integral of the average kernel for Brownian coagulation in the continuous regime is given by

$$II_{CR} = \int_0^\infty \int_0^\infty \left(\eta_1^{-\frac{1}{3}} + \eta^{-\frac{1}{3}}\right)\left(\eta_1^{\frac{1}{3}} + \eta^{\frac{1}{3}}\right) \exp(-\eta - \eta_1)\, d\eta\, d\eta_1 \tag{36}$$

By expanding the homogeneous power-law kernel directly, the double integral can be reorganized into the second Euler integral as

$$II_{CR} = \int_0^\infty \int_0^\infty \left(2 + \eta^{\frac{1}{3}}\eta_1^{-\frac{1}{3}} + \eta^{-\frac{1}{3}}\eta_1^{\frac{1}{3}}\right) \exp(-\eta - \eta_1)\, d\eta\, d\eta_1 \tag{37}$$

Using the definition of the Euler gamma function,

$$\Gamma(z) = \int_0^\infty t^{z-1} \exp(-t)\, dt \tag{38}$$

Then, the exact expression of the average kernel for Brownian coagulation in the continuous regime is:

$$II_{CR} = 2\Gamma(1) + 2\Gamma\left(\frac{2}{3}\right)\Gamma\left(\frac{4}{3}\right) \tag{39}$$

The result of the double integrals with nine significant digits was $II_{CR} = 4.41839915$.

**Shear coagulation**

The double integral of the average kernel for shear coagulation is given by

$$II_{SC} = \int_0^\infty \int_0^\infty \left(\eta_1^{\frac{1}{3}} + \eta^{\frac{1}{3}}\right)^3 \exp(-\eta - \eta_1)\, d\eta\, d\eta_1 \tag{40}$$

Analogously, by expanding the homogeneous power-law kernel directly, the double integral can be reorganized into the second Euler integral as

$$II_{SC} = \int_0^\infty \int_0^\infty \left(\eta + \eta_1 + \eta^{\frac{1}{3}}\eta_1^{-\frac{1}{3}} + \eta^{-\frac{1}{3}}\eta_1^{\frac{1}{3}}\right) \exp(-\eta - \eta_1)\, d\eta\, d\eta_1 \tag{41}$$

And the corresponding exact expression of the average kernel for shear coagulation is

$$II_{SC} = 2\Gamma(2) + 6\Gamma\left(\frac{5}{3}\right) \cdot \Gamma\left(\frac{4}{3}\right) \tag{42}$$

The result of the double integrals with nine significant digits was $II_{SC} = 6.83679830$.

**Brownian coagulation in the free molecule regime**

The double integral of the average kernel for Brownian coagulation in the free molecular regime is

$$II_{FM} = \int_0^\infty \int_0^\infty \left(\frac{1}{\eta} + \frac{1}{\eta_1}\right)^{1/2} \left(\eta^{\frac{1}{3}} + \eta_1^{\frac{1}{3}}\right)^2 \exp(-\eta - \eta_1)\, d\eta\, d\eta_1 \tag{43}$$

Let $\eta = a^2$ and $\eta_1 = b^2$, then, $d\eta = 2a\,da$ and $d\eta_1 = 2b\,db$. Substituting them into the double integrals, it becomes

$$II_{FM} = 4\int_0^\infty \int_0^\infty \sqrt{a^2 + b^2}\left(a^{\frac{2}{3}} + b^{\frac{2}{3}}\right)^2 \exp(-a^2 - b^2)\, da\, db \tag{44}$$

Let $a = r\cos\theta, b = r\sin\theta, (r \geq 0, \theta \in [0, \pi/2])$, then, the integral above is transformed to

$$II_{FM} = 4\int_0^{\frac{\pi}{2}} \int_0^\infty r^{\frac{10}{3}} \exp(-r^2)(\cos^{\frac{2}{3}}\theta + \sin^{\frac{2}{3}}\theta)^2\, dr\, d\theta \tag{45}$$

The integral can be expressed as the product of the radial integrals $I_r$ and the angular integrals $I_\theta$. For the radial integral, it can be expressed directly as

$$I_r = \int_0^\infty r^{\frac{10}{3}} \exp(-r^2)\, dr = \frac{1}{2}\Gamma\left(\frac{13}{6}\right) \tag{46}$$

For the angular integral, it can be expressed as

$$I_\theta = \int_0^{\frac{\pi}{2}}(\cos^{\frac{2}{3}}\theta + \sin^{\frac{2}{3}}\theta)^2\, d\theta = 2B\left(\frac{7}{6}, \frac{1}{2}\right) + 2B\left(\frac{5}{6}, \frac{5}{6}\right) \tag{47}$$

where $B$ denotes the Beta function. According to the property of the Beta function

$$B(x, y) = \frac{\Gamma(x)\Gamma(y)}{\Gamma(x+y)} \tag{48}$$

Then the double integral of the average kernel for Brownian coagulation in the free molecule regime can be rewritten as a combination of Gamma functions as

$$II_{FM} = \Gamma\left(\frac{13}{6}\right)\frac{\Gamma\left(\frac{7}{6}\right)\Gamma\left(\frac{1}{2}\right) + \Gamma\left(\frac{5}{6}\right)\Gamma\left(\frac{5}{6}\right)}{\Gamma\left(\frac{10}{6}\right)} \tag{49}$$

The result of the double integrals with nine significant digits was $II_{FM} = 6.99822397$.

**Sedimentary coagulation**

The double integral of the average kernel for the sedimentary coagulation is

$$II_{SD} = \int_0^\infty \int_0^\infty (\eta^{\frac{1}{3}} + \eta_1^{\frac{1}{3}})^3 \left|\eta^{\frac{1}{3}} - \eta_1^{\frac{1}{3}}\right| \exp(-\eta - \eta_1)\, d\eta\, d\eta_1 \tag{50}$$

Let $x = \eta^{\frac{1}{3}}$ and $y = \eta_1^{\frac{1}{3}}$, $d\eta = 3x^2 dx$ and $d\eta_1 = 3y^2 dy$. Substituting them into the double integral, it becomes

$$II_{SD} = \int_0^\infty \int_0^\infty 9x^2 y^2 (x + y)^3 |x - y| \exp(-x^3 - y^3)\, dx\, dy \tag{51}$$

The double integral above can be divided into two parts as

$$\begin{cases} II_1 = \int_0^\infty \int_0^x 9(x+y)^3(x-y)x^2y^2 \exp(-x^3-y^3)\,dxdy & (x>y) \\ II_2 = \int_0^\infty \int_0^y 9(x+y)^3(y-x)x^2y^2 \exp(-x^3-y^3)\,dydx & (y>x) \end{cases} \quad (52)$$

It is evident that $II_1 = II_2$, thus, $II_{SD} = 2II_1$. Applying the polar coordinate transformation ($x = r\cos\theta$ and $y = r\sin\theta$) and the definition of the generalized gamma function

$$\Gamma(\alpha,\lambda,\mu) = \int_0^\infty x^\alpha \exp(-\lambda x^\mu) = \left(\mu \lambda^{\frac{\alpha+1}{\mu}}\right)^{-1} \Gamma\left(\frac{\alpha+1}{\mu}\right) \quad (53)$$

The double integral $II_{SD}$ can be simplified as

$$II_{SD} = 6\Gamma\left(\frac{10}{3}\right) \int_0^{\frac{\pi}{4}} \frac{(\cos\theta+\sin\theta)^3(\cos\theta-\sin\theta)\cos^2\theta\sin^2\theta}{(\cos^3\theta+\sin^3\theta)^{\frac{10}{3}}} d\theta \quad (54)$$

Using the trigonometric substitution $t = \tan\theta$, the exact expression of the average kernel for sedimentary coagulation becomes

$$II_{SD} = 6\Gamma\left(\frac{10}{3}\right) \int_0^1 \frac{t^2(1-t)}{(1-t+t^2)^3(1+t^3)^{\frac{1}{3}}} dt \quad (55)$$

The result of the double integrals with nine significant digits was $II_{SD} = 2.58940496$.

**The approximation of double integrals with Gauss-Laguerre quadrature**

Based on the MATLAB code of the Gauss-Laguerre quadrature in **Appendix B**, the corresponding values of the double integrals $II$, together with the error $\epsilon_{n,n}$, constant slope $C$, $Q_{n+1,n+1}$ and its remainder $R_{n+1,n+1}$ are listed in Table 1. The change in the value $Q_{n+1,n+1}$ and its error $\epsilon_{n,n}$ as the number of Laguerre polynomial nodes increases are shown in Fig.1. The results show that the simplified error analysis of double integrals with the Gauss-Laguerre quadrature in this study is reasonable.

The value of the pre-exponential factor $p$ of the average kernel was half that of the double integrals $II$. The results of the average kernel values calculated using the two-dimensional Gaussian-Laguerre (GL) quadrature method and direct computation (DC) for different collision kernels are compared in Table 2. The results demonstrate that the computational method of the average kernel in this study is reliable and universal.

Table 1. The numerical values related to the Gauss-Laguerre quadrature for double integrals

| Type | $\epsilon_{n,n}$ | $C$ | $Q_{n+1,n+1}$ | $R_{n+1,n+1}$ | $II$ |
|---|---|---|---|---|---|
| CR | $2.5865E-5$ | $-1.6572$ | $4.4025$ | $0.0141$ | $4.4025 \pm 0.0141$ |
| SC | $9.3561E-7$ | $-2.3733$ | $6.8371$ | $0.0002$ | $6.8371 \pm 0.0002$ |
| FM | $1.2344E-4$ | $-1.5209$ | $6.9032$ | $0.0852$ | $6.9032 \pm 0.0852$ |
| SD | $7.4036E-6$ | $-1.9750$ | $2.5861$ | $0.0027$ | $2.5861 \pm 0.0027$ |

Table 2. The converged average kernels for a few common collision kernels

| Type | $\beta(v,v_1)$ | $\bar{\beta} = pu^q$ (GL) | $\bar{\beta} = pu^q$ (DC) |
|---|---|---|---|
| CR | $\beta_{CR} = \left(v_1^{-\frac{1}{3}} + v^{-\frac{1}{3}}\right)\left(v_1^{\frac{1}{3}} + v^{\frac{1}{3}}\right)$ | $\bar{\beta}_{CR} = 2.2013$ | $\bar{\beta}_{CR} = 2.2092$ |
| SC | $\beta_{SC} = \left(v_1^{\frac{1}{3}} + v^{\frac{1}{3}}\right)^3$ | $\bar{\beta}_{SC} = 3.4186u$ | $\bar{\bar{\beta}}_{SC} = 3.4184u$ |
| FM | $\beta_{FM} = (v_1^{-1} + v^{-1})^{\frac{1}{2}}\left(v_1^{\frac{1}{3}} + v^{\frac{1}{3}}\right)^2$ | $\bar{\beta}_{FM} = 3.4516u^{\frac{1}{6}}$ | $\bar{\bar{\beta}}_{FM} = 3.4991u^{\frac{1}{6}}$ |
| SD | $\beta_{SD} = \left(v_1^{\frac{1}{3}} + v^{\frac{1}{3}}\right)^3 \left|v_1^{\frac{1}{3}} - v^{\frac{1}{3}}\right|$ | $\bar{\beta}_{SD} = 1.2931u^{\frac{4}{3}}$ | $\bar{\bar{\beta}}_{SD} = 1.2947u^{\frac{4}{3}}$ |

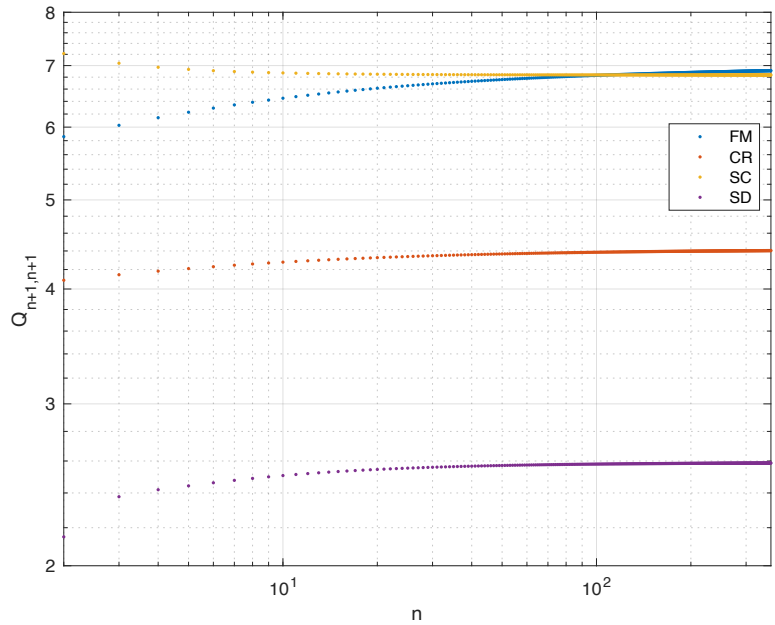

a

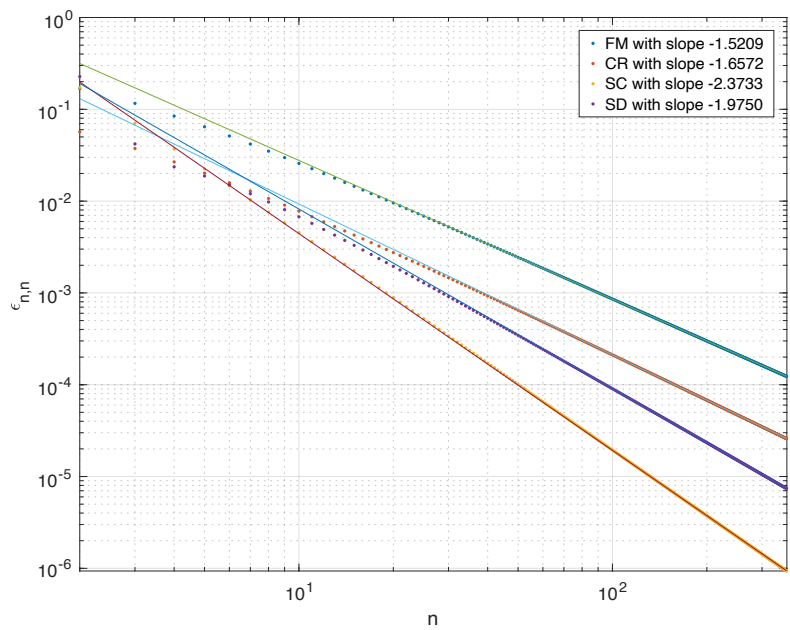

b

Fig.1. Numerical results of Gauss-Laguerre quadrature for double integral. a) the change in the value $Q_{n+1,n+1}$ with the number of interpolation nodes $n$; b) the change in the error $\epsilon_{n,n}$ as the number of Laguerre polynomial nodes increases.

**Conclusions**

In this study, a compact and fast MATLAB code based on the Gauss-Laguerre quadrature for double integrals is proposed to calculate the pre-exponential factor of the average kernel after proper normalization, which reduces the complexity of the traditional average kernel method in the literature. The proposed numerical algorithm provides a robust computational basis for extending and applying the AK-iDNS framework to solve the Smoluchowski equation under various conditions.

**Acknowledgement**

This research was funded by the National Natural Science Foundation of China (Grant number 11972169).

**Nomenclature**

| | |
|---|---|
| $x$ | horizontal axis in the Cartesian coordinate |
| $y$ | vertical axis in the Cartesian coordinate |
| $z$ | complex variable |
| $r$ | radial axis in the polar coordinates |
| $\theta$ | circumferential axis in the polar coordinates |
| $n$ | the number density of particles |
| $v$ | particle volume |
| $u$ | the algebraic mean volume |
| $t$ | time or temporary variable |
| $\beta$ | the collision kernel of coagulation |
| $\bar{\beta}$ | the average collision kernel |
| $f$ | the function |
| $m, n, r, s$ | the order parameters |
| $\alpha$ | the scaler factor or power index |
| $\alpha_0$ | the constant |
| $p$ | the proportional factor |
| $q$ | the power index |
| $I$ | single integrals |
| $II$ | single integrals |
| $L$ | Laguerre polynomials |
| $Q_n$ | approximate value of single integrals |
| $Q_{n,n}$ | approximate value of double integrals |
| $R_n$ | remainder of single integrals |
| $R_{n,n}$ | remainder of double integrals |
| $\epsilon_n$ | the error of single integrals |
| $\epsilon_{n,n}$ | the error of double integrals |
| $\Gamma$ | Euler gamma function |
| $B$ | Beta function |
| $\mu, \lambda$ | the parameter in Euler gamma function |
| $M_0$ | total particle number concentration |
| $M_1$ | the particle volume concentration |
| $\eta$ | the dimensionless particle volume, |
| $\psi$ | the dimensionless particle size distribution |
| $C$ | the constant |
| $C_z$ | contour integrals |
| $U$ | the second kind confluent hypergeometric function |
| $k$ | function |
| $\Pi$ | function |
| AK | average kernel method |
| iDNS | iterative direct numerical simulation |
| CR | Brownian coagulation in the continuum regime |
| FM | Brownian coagulation in the free molecule regime |
| SC | shear coagulation |
| SD | sedimentary coagulation |
| GL | Gauss-Laguerre quadrature |
| DC | direct computation |

**Appendix A: Population average kernel**

In the literature, the population-averaged coagulation kernel is defined as

$$\int_0^\infty \int_0^\infty \bar{\beta}\, n(v_1)n(v)\, dv_1 dv = \int_0^\infty \int_0^\infty \beta(v_1, v)\, n(v_1)n(v)\, dv_1 dv \quad (A1)$$

It can be calculated as

$$\bar{\beta} = \frac{\int_0^\infty \int_0^\infty \beta(v_1, v)\, n(v_1)n(v)\, dv_1 dv}{\int_0^\infty \int_0^\infty n(v_1)n(v)\, dv_1 dv} \quad (A2)$$

Based on the analytical similarity solution (**Schumann, 1940; Pan et al., 2024**),

$$\psi(\eta) = e^{-\eta} \quad (A3)$$

The particle number density can be expressed as

$$n(v, t) = \frac{N}{u} \exp\left[-\frac{v}{u}\right] \quad (A4)$$

where $N = \int_0^\infty n(v, t) dv$ is the total particle number. And the population-averaged kernel can be expressed as

$$\bar{\beta} = \frac{1}{u^2} \int_0^\infty \int_0^\infty \beta(v_1, v) \exp\left(-\frac{v_1+v}{u}\right) dv_1 dv \quad (A5)$$

Therefore, the average kernel based on Laplace transformation is equivalent to the population-averaged kernel.

**Appendix B: Evaluation points and weights of the Gauss-Laguerre quadrature**

The Laguerre polynomials were calculated recursively. Defining the first two polynomials as

$$L_0(x) = 1, \ L_1(x) = 1 - x; \tag{B1}$$

And the recurrence relation of Laguerre polynomials for any $n \geq 1$ is,

$$L_{n+1}(x) = \frac{(2n+1-x)L_n(x) - nL_{n-1}(x)}{n+1}; \tag{B2}$$

It can be written explicitly as

$$L_{n+1}(x) = \sum_{k=0}^{n+1} (-1)^k \frac{n!}{(k!)^2 (n-k)!} x^k \tag{B3}$$

The abscissas $x_i$ can be obtained from the zeros of the Laguerre polynomials, i.e., $L_{n+1}(x) = 0$. The corresponding weights $w_i$ can be calculated explicitly as

$$w_i = \frac{x_i}{(n+1)^2 [L_{n+1}(x_i)]^2} \tag{B4}$$

The corresponding MATLAB code is as follows.

```
function [x,w] = Gauss_Laguerre(n)
% n: Define the order of the Gauss-Laguerre quadrature
% x: the evaluation points of the nth Laguerre polynomial
% w: the weights of the nth Gauss-Laguerre quadrature
% w_i = x_i /((n+1)^2*(L_{n+1}(x_i))^2)
syms y;
p1 = laguerreL(n, y);
x = double(solve(p1 == 0, y));
w = zeros(n, 1);
for i = 1:n
    L_n_plus_1 = laguerreL(n+1, x(i));
    w(i) = double((x(i)/((n+1)^2*(L_n_plus_1)^2)));
end
end
```

For instance, the Laguerre polynomials $L_{n+1}(x)$ for $n = 10$ is

$$L_{11}(x) = \frac{x^{10}}{3628800} - \frac{x^9}{36288} + \frac{x^8}{896} - \frac{x^7}{42} +$$

$$\frac{7x^6}{24} - \frac{21x^5}{10} + \frac{35x^4}{4} - 20x^3 + \frac{45x^2}{2} - 10x + 1 \tag{B5}$$

The corresponding evaluation points and weights for the Gauss–Lagerre quadrature are presented in Table A1.

Table A1. The evaluation points and weights of Gauss Laguerre quadrature for $n = 10$

| $i$ | $x_i$ | $A_i$ |
|---|---|---|
| 1 | 0.1377 | 0.3084 |
| 2 | 0.7294 | 0.4011 |
| 3 | 1.8083 | 0.2180 |
| 4 | 3.4014 | 0.0620 |
| 5 | 5.5524 | 0.0095 |
| 6 | 8.3301 | 0.0007 |
| 7 | 11.8437 | 2.8E-05 |
| 8 | 16.2792 | 4.2E-07 |
| 9 | 21.9965 | 1.8E-09 |
| 10 | 29.9206 | 9.9E-13 |

Taking the shear coagulation kernel as an example, the MATLAB code for the calculation of $Q_{n+1,n+1}$ and its error $\epsilon_{n,n}$ based on the Gauss-Laguerre quadrature for the double integrals is as follows:

```
clear,
format long
n = 99;
for k = 1:n
    beta = zeros(k,k); [x,w] = Gauss_Laguerre(k); Q = 0;
    for i = 1:k
        for j = 1:k
            % collision kernel for shear coagulation
            beta(i,j) = ( x(i)^(1/3) + x(j)^(1/3) )^3;
            Q = Q + w(i)*w(j)*beta(i,j);
        end
    end
    p(k) = Q/2;
end
error = abs(p(2:n)-p(1:n-1));
figure, loglog(1:n,p,'.'),grid on
figure, loglog(1:n-1,error,'.'),grid on
end
```